\documentstyle[amsmath,amsfonts,amssymb,amsthm,amscd,12pt]{article}
\font\sixbb=msbm6
\font\eightbb=msbm8
\font\twelvebb=msbm10 scaled 1095
\newfam\bbfam
\textfont\bbfam=\twelvebb \scriptfont\bbfam=\eightbb
                           \scriptscriptfont\bbfam=\sixbb
\def\bb{\fam\bbfam\twelvebb}

\newcommand{\Com}{{\bb C}}
\newcommand{\Int}{{\bb Z}}

\newcommand{\Rat}{{\bb Q}}

\newtheorem{theorem}{\bf Theorem}

\newtheorem{proposition}[theorem]{\bf Proposition}

\newcommand{\enp}{\begin{flushright} $\Box$ \end{flushright}}
\newcommand{\beq}[0]{\begin{equation}}
\newcommand{\enq}[0]{\end{equation}}
\newcommand{\rk}{rank}

\newcommand{\Ght}{\widehat{G}}
\newcommand{\fht}{\widehat{f}}
\newcommand{\hht}{\widehat{h}}
\newcommand{\cf}{{\cal F}}
\newcommand{\cb}{{\cal B}}
\newcommand{\spp}{{\rm supp}}
\newcommand{\hh}{{\rm H}}
\newcommand{\thh}{\tilde{\rm H}}
\newcommand{\hb}{{\rm B}}

\newcommand{\one}{{\bf 1}}
\newcommand{\ch}{{\cal H}}

\title{Homology of Balanced Complexes  \\ via the Fourier Transform}
\begin{document}
\author{Roy Meshulam\thanks{Department of Mathematics,
Technion, Haifa 32000, Israel. e-mail:
meshulam@math.technion.ac.il~. Supported by ISF and BSF grants. }
}
\maketitle
\pagestyle{plain}
\begin{abstract}
Let $G_0,\ldots,G_k$ be finite abelian groups and let $G_0* \cdots *G_k$ be the join of the $0$-dimensional complexes $G_i$. We give a characterization of the integral $k$-coboundaries of subcomplexes of $G_0* \cdots *G_k$ in terms of the Fourier transform on the group $G_0 \times \cdots \times G_k$. This leads to an extension of a recent result of Musiker and Reiner  on a topological interpretation of the cyclotomic polynomial.
\end{abstract}

\section{Introduction}
\ \ \ \
Let $G_0,\ldots,G_k$ be finite abelian groups with the discrete topology and let $N=\prod_{i=0}^k (|G_i|-1)$. The simplicial join $Y=G_0* \cdots *G_k$ is homotopy equivalent to a wedge of $N$ $k$-dimensional spheres. Subcomplexes of
$Y$ are called {\it balanced complexes} (see e.g. \cite{Stanley}). Denote the $(k-1)$-dimensional skeleton of $Y$ by $Y^{(k-1)}$.
Let $A$ be a subset of $G_0 \times \cdots \times G_k$.  Regarding each $a \in A$ as an oriented $k$-simplex
of $Y$, we consider the balanced complex
$$X(A)=X_{G_0,\ldots,G_k}(A)=Y^{(k-1)} \cup A.$$
In this note we characterize the integral $k$-coboundaries of $X(A)$ in terms of the Fourier transform on the group $G_0 \times \cdots \times G_k$. As an application we give a short proof of an extension of a recent result of Musiker and Reiner \cite{MR10} on a topological interpretation of the cyclotomic polynomial.

We recall some terminology. Let $R[G]$ denote the group algebra of a finite abelian group $G$ with coefficients in a ring $R$. By writing  $f=\sum_{x \in G} f(x) x\in R[G]$ we identify elements of $R[G]$ with $R$-valued functions on $G$.
For a subset $A \subset G$ let $R[A]=\{f \in R[G]: \spp(f) \subset A\}$.
Let $\Ght$ be the character group of $G$.
The Fourier transform is the linear bijection $\cf:\Com[G] \rightarrow \Com[\Ght]$
given on $f \in \Com[G]$ and $\chi \in \Ght$ by
$$ \cf(f)(\chi)=\fht(\chi)=\sum_{x \in G}
f(x) \chi(x)~.
$$
Let $G= G_0 \times \cdots \times G_k$ then $\Ght=\Ght_0 \times \cdots \times \Ght_k$.
For $0 \leq i \leq k$ let
$$L_i=G_0 \times \cdots \times G_{i-1} \times G_{i+1} \times \cdots \times G_k.$$
We identify the group of integral $k$-cochains $C^k(X(A);\Int)$ with $\Int[A]$ and the group of integral $(k-1)$-cochains $C^{k-1}(X(A);\Int)=C^{k-1}(X(G);\Int)$ with the $(k+1)$-tuples
$\psi=(\psi_0,\ldots,\psi_k)$ where $\psi_i \in \Int[L_i]$. The coboundary map
$$d_{k-1}:C^{k-1}(X(G);\Int) \rightarrow C^k(X(G);\Int)$$ is given by
$$d_{k-1} \psi(g_0,\ldots,g_k)=\sum_{i=0}^k (-1)^i \psi_i(g_0,\ldots,g_{i-1},g_{i+1},\ldots,g_k).$$
For $0 \leq i \leq k$ let $\one_i$ denote the trivial character of $G_i$ and
let $$\Ght^+=(\Ght_0-\{\one_0\}) \times \cdots \times (\Ght_k-\{\one_k\}) .$$
For $A \subset G$ and $f \in \Int[G]$ let $f_{|A} \in \Int[A]$ denote the restriction of $f$ to $A$.
The group $\hb^k(X(A);\Int)=d_{k-1}C^{k-1}(X(G);\Int)$ of integral $k$-coboundaries of $X(A)$ is characterized by the following
\begin{proposition}
\label{cob}
For any $A \subset G$
$$\hb^k(X(A);\Int) = \{f_{|A}: f \in \Int[G] ~~such~that~~ \spp(\fht) \subset \Ght -\Ght^+\}.$$
\end{proposition}
As an application of Proposition \ref{cob} we study the homology of the following family of balanced complexes introduced by Musiker and Reiner \cite{MR10}.
Let $p_0,\ldots,p_k$ be distinct primes and for $0 \leq i \leq k$ let $G_i=\Int/p_i \Int=\Int_{p_i}$.
Writing $n=\prod_{i=0}^k p_i$, let $\theta:\Int_n \rightarrow G=G_0 \times \cdots \times  G_k$ be the  standard isomorphism given by
$$\theta(x)=(x({\rm mod}~p_0),\ldots, x({\rm mod}~p_k)).$$
For any $\ell$ let $\Int_{\ell}^*=\{m \in \Int_{\ell}: \text{gcd}(m,{\ell})=1\}$.
Let $\varphi(n)=|\Int_n^*|=\prod_{i=0}^k (p_i-1)$ be the Euler function of $n$ and
let $A_0=\{\varphi(n)+1,\ldots,n-1\}$. For $A \subset \{0,\ldots,\varphi(n)\}$ consider the complex
$$K_A=X(\theta(A \cup A_0)) \subset \Int_{p_0} * \cdots * \Int_{p_k}~.$$

Let $\omega=\exp(\frac{2\pi i}{n})$ be a fixed primitive $n$-th root of unity. The $n$-th cyclotomic polynomial
(see e.g. \cite{Lang}) is given by $$\Phi_n(z)=\prod_{j \in \Int_n^*} (z-w^j)=\sum_{j=0}^{\varphi(n)} c_j z^j \in \Int[z].$$
Musiker and Reiner \cite{MR10} discovered the following remarkable connection between the coefficients of $\Phi_n(z)$ and the homology of the complexes $K_{\{j\}}$.
\begin{theorem}[Musiker and Reiner]
\label{mr}
For any $j \in  \{0,\ldots,\varphi(n)\}$
$$\thh_i(K_{\{j\}};\Int) \cong
\left\{
\begin{array}{ll}
\Int/c_j \Int  & i=k-1 \\
\Int & i=k~\text{and}~c_j=0 \\
0 & \text{otherwise.}
\end{array}
\right.~~
$$
\end{theorem}
\noindent
The next result extends Theorem \ref{mr}
to general $K_A$'s. Let $$c_A=(c_j: j \in A) \in \Int^{A}$$ and
$$d_A=
\left\{
\begin{array}{ll}
\text{gcd}(c_A)  & c_A \neq 0 \\
0 & c_A=0
\end{array}
\right.~~ $$

\begin{theorem}
\label{hka}
For any $A \subset \{0,\ldots,\varphi(n)\}$
$$
\thh^i(K_A;\Int)\cong
\left\{
\begin{array}{ll}
\Int & i=k-1 ~\text{and}~ d_A=0 \\
\Int^{|A|-1}\oplus \Int/d_A \Int & i=k \\
0 & \text{otherwise}
\end{array}
\right.~~
$$
and
$$
\thh_i(K_A;\Int)\cong
\left\{
\begin{array}{ll}
\Int/d_A\Int & i=k-1 \\
\Int^{|A|} & i=k ~\text{and}~ d_A=0 \\
\Int^{|A|-1} & i=k ~\text{and}~ d_A \neq 0 \\
0 & \text{otherwise.}
\end{array}
\right.~~
$$
\end{theorem}
\ \\ \\

Proposition \ref{cob} is proved in Section \ref{s:cob}. It is then used in Section \ref{s:hka} to obtain an explicit
form of the $k$-coboundaries of $K_A$ (Proposition \ref{isoh}) that directly implies Theorem \ref{hka}.

\section{$k$-Coboundaries and Fourier Transform}
\label{s:cob}
\noindent
{\bf Proof of Proposition \ref{cob}:} It suffices to consider the case $A=G$.
Let $\psi=(\psi_0,\ldots,\psi_k) \in C^{k-1}(X(G);\Int)$.
Then for any $\chi=(\chi_0, \ldots,\chi_k) \in \Ght$
$$\widehat{d_{k-1} \psi}(\chi)=\sum_{g=(g_0,\ldots,g_k) \in G} d_{k-1} \psi(g) \chi(g)=$$
$$\sum_{(g_0,\ldots,g_k)} \sum_{i=0}^k (-1)^i \psi_i(g_0,\ldots,g_{i-1},g_{i+1},\ldots,g_k) \prod_{j=0}^k \chi_j(g_j)=$$
$$\sum_{i=0}^k (-1)^i \sum_{(g_0,\ldots,g_{i-1},g_{i+1},\ldots ,g_k)}\psi_i(g_0,\ldots,g_{i-1},g_{i+1},\ldots,g_k)
\prod_{j \neq i} \chi_j(g_j)
\sum_{g_i} \chi_i(g_i)=$$
$$\sum_{i=0}^k (-1)^i \widehat{\psi_i}(\chi_0,\ldots,\chi_{i-1},\chi_{i+1},\ldots,\chi_k) |G_i| \delta(\chi_i,\one_i)$$
where $\delta(\chi_i,\one_i)=1$ if $\chi_i=\one_i$ and is zero otherwise.
\\
Therefore $\spp(\widehat{d_{k-1} \psi}) \subset \Ght -\Ght^+$ and so
$$U_1\stackrel{\text{def}}{=}\hb^k(X(G);\Int) \subset \{f \in \Int[G]:\spp(\fht) \subset \Ght -\Ght^+\}\stackrel{\text{def}}{=}U_2.$$

Since $X(G)$ is homotopy equivalent to a wedge of $\prod_{i=0}^k (|G_i|-1) =|\Ght^+|~$
$k$-dimensional spheres, it follows  $\hh^k(X(G);\Int)=\Int[G]/U_1$ is free of rank $|\Ght^+|$ and hence
$\rk~U_1=|\Ght|-|\Ght^+|$. On the other hand, the injectivity of the Fourier transform implies that
$$\rk~ U_2 \leq \dim_{\Com} \{f \in \Com[G]:\spp(\fht) \subset \Ght -\Ght^+\}=|\Ght|-|\Ght^+|$$
and therefore  $\rk~U_2/U_1=0$. Since $U_2/U_1 \subset \hh^k(X(G);\Int)$ is free it follows that $U_1=U_2$.
{\enp}

\section{The Homology of $K_A$}
\label{s:hka}
\ \ \ \
Recall that $G=\Int_{p_0} \times \cdots \times \Int_{p_k}$ and $n=\prod_{j=0}^k p_j$.
For $h \in \Int[G]$  let $\theta^{*}h \in \Int[\Int_n]$ be the pullback of $h$ given by
$\theta^*h(x)=h(\theta(x))$.
For any $\ell$ we identify the character group $\widehat{\Int_{\ell}}$ with $\Int_{\ell}$ via the isomorphism
$\eta_{\ell}:\Int_{\ell} \rightarrow \widehat{\Int_{\ell}}$ given by $\eta_{\ell}(x)(y)= \exp(2 \pi i x y/{\ell}).$

Proposition \ref{cob} implies the following characterization of the integral
 $k$-coboundaries of $K_A$.
For $A\subset \{0,\ldots,\varphi(n)\}$ let
$\theta_A$ denote the restriction of $\theta$ to $A \cup A_0$ and let $\theta_A^*$ be the induce isomorphism from $\Int[\theta(A \cup A_0)]$ to $\Int[A \cup A_0].$
Let
$$\cb(A)=\{f_{|A \cup A_0}: f \in \Int[\Int_n] ~~{such~that}~~
\fht(1)=0 \}.$$
\begin{proposition}
\label{isoh}
$$
\theta_A^*\hb^k(K_A;\Int)=\cb(A).
$$
\end{proposition}
\noindent
{\bf Proof:} We first examine the relation between the Fourier transforms on $\Int_n$ and  on
$G$.
Let $$\lambda=\sum_{j=0}^k \prod_{t \neq j} p_t \in \Int_n^*.$$
For any $h \in \Int[G]$ and $m \in \Int_n$
$$\widehat{\theta^*h}(\lambda m)=\sum_{x \in \Int_n} \theta^*h(x) \exp(\frac{2 \pi i x \lambda m}{n})=
$$
\begin{equation}
\label{fourh}
\sum_{x \in \Int_n} h(\theta(x))\exp(\sum_{j=0}^k \frac{2 \pi i x m}{p_j})=\hht(\theta(m)).
\end{equation}
\noindent
Noting that  $$\theta^{-1}(\Ght^+)=\theta^{-1}(\Int_{p_0}^* \times \cdots \times \Int_{p_k}^*)=\Int_n^*=\lambda \Int_n^*$$ it follows from
Proposition \ref{cob} and Eq. (\ref{fourh}) that
$$\hb^k(K_A;\Int) =\{h_{|\theta(A \cup A_0)}: h \in \Int[G] ~~{such~that}~~\spp(\hht) \subset \Ght-\Ght^+\}=
$$
\begin{equation}
\label{tran}
(\theta_A^*)^{-1}\{f_{|A \cup A_0}: f \in \Int[\Int_n] ~~{such~that}~~\spp(\fht) \subset \Int_n - \Int_n^* \}.
\end{equation}
\ \\ \\
Let ${\cal P}_n=\{\omega^m: m \in \Int_n^*\}$ be the set of primitive $n$-th roots of $1$.
The Galois group ${\rm Gal}(\Rat(\omega)/\Rat)$ acts transitively
on ${\cal P}_n$.  Hence, by Eq. (\ref{tran}):
$$
\theta_A^*\hb^k(K_A;\Int) = \{f_{|A \cup A_0}: f \in \Int[\Int_n] ~~{such~that}~~\spp(\fht) \subset \Int_n - \Int_n^* \}=
$$
$$
\{f_{|A \cup A_0}: f \in \Int[\Int_n] ~~{such~that}~~
\fht(m)=0 {~~for~all~~} m \in \Int_n^*\}=
$$
$$
\{f_{|A \cup A_0}: f \in \Int[\Int_n] ~~{such~that}~~
\fht(1)=0 \}=\cb(A).
$$
{\enp}
\ \\ \\
{\bf Proof of Theorem \ref{hka}:} Proposition \ref{isoh} implies that $\theta_A^*$ induces an isomorphism between $\hh^k(K_A;\Int)$ and
$$\ch(A) \stackrel{\text{def}}{=} \Int[A \cup A_0]/\cb(A).$$
For $j \in A \cup A_0$ let $g_j \in \Int[A \cup A_0]$ be given by
$g_j(i)=1$ if $i=j$ and $g_j(i)=0$ otherwise. Let $[g_j]$ be the image of $g_j$ in $\ch(A)$.
The computation of $\ch(A)$ depends on the following two observations:
\ \\ \\
(i) $\ch(A)$ is generated by $\{[g_j]: j \in A\}$.
\ \\ \\
{\bf Proof:} Let $t \in A_0$. There exist $u_0,\ldots,u_{\varphi(n)-1} \in \Int$ such that
$$\sum_{\ell=0}^{\varphi(n)-1} u_{\ell} \omega^{\ell} +\omega^t=0.$$
Let $f \in \Int[\Int_n]$ be given by
$$f(\ell)=
\left\{
\begin{array}{ll}
u_{\ell} & 0 \leq \ell \leq \varphi(n)-1 \\
1 & \ell=t \\
0 & \text{otherwise.}
\end{array}
\right.~~
$$
Since
$$\fht(1)=\sum_{\ell=0}^{\varphi(n)-1} u_{\ell} \omega^{\ell} +\omega^t=0$$
it follows that
$$
\sum_{j \in A} u_j g_j+g_t=f_{|A \cup A_0} \in \cb(A).
$$
Hence  $[g_t]=-\sum_{j \in A} u_j [g_j]$.
\ \\ \\
(ii) The minimal relation between $\{[g_j]\}_{j \in A}$ is $\sum_{j \in A} c_j[g_j]=0$.
\ \\ \\
{\bf Proof:} Let $f \in \Int[\Int_n]$ be given by
$f(\ell)=c_{\ell}$ if $0 \leq \ell \leq \varphi(n)$ and zero otherwise.
Since $\fht(1)=\Phi_n(\omega)=0$, it follows that
$$
\sum_{j \in A}c_j g_j = f_{|A \cup A_0} \in \cb(A).
$$
Hence $\sum_{j \in A} c_j [g_j]=0$.
On the otherhand, if $\sum_{j \in A} \alpha_j [g_j]=0$ then there exists an $h \in \Int[\Int_n]$
such that $\hht(1)=0$ and $h_{|A \cup A_0}=\sum_{j \in A} \alpha_j g_j$. In particular $h(\ell)=0$ for $\ell \geq \varphi(n)+1$.
Let $p(z)=\sum_{\ell=0}^{\varphi(n)} h(\ell) z^{\ell}$  then $p(\omega)=\hht(1)=0$.
Hence  $p(z)=r\Phi_n(z)$ for some $r \in \Int$. Therefore
$\alpha_j=h(j)=rc_j$ for all $j \in A$.
\ \\ \\
It follows from (i) and (ii) that
\begin{equation}
\label{cohomk}
\hh^k(K_A;\Int)\cong \ch(A)=\Int[A]/\Int c_A
\cong \Int^{|A|-1}\oplus \Int/d_A \Int ~.
\end{equation}
The remaining parts of Theorem \ref{hka} are formal consequences of (\ref{cohomk}) and the universal coefficient theorem (see e.g. \cite{Munkres}):
\begin{equation}
\label{uct}
0 \leftarrow {\rm Hom}(\hh_p(K_A;\Int),\Int) \leftarrow \hh^p(K_A;\Int) \leftarrow {\rm Ext}(\hh_{p-1}(K_A;\Int),\Int) \leftarrow 0~.
\end{equation}
First consider the case $c_A=0$. By (\ref{cohomk}) and (\ref{uct})
$$
0 \leftarrow {\rm Hom}(\hh_k(K_A;\Int),\Int) \leftarrow \Int^{|A|} \leftarrow {\rm Ext}(\hh_{k-1}(K_A;\Int),\Int) \leftarrow 0~.
$$
Therefore $\hh_k(K_A;\Int)\cong \Int^{|A|}$ and $\hh_{k-1}(K_A;\Int)$ is torsion free. The
Euler-Poincar\'{e} relation
\begin{equation}
\label{EPH}
\rk~ \hh_k(K_A;\Int)=\rk~ \thh_{k-1}(K_A;\Int)+|A|-1~
\end{equation}
then implies
that $\thh_{k-1}(K_A;\Int)\cong \Int$ and $$\thh^{k-1}(K_A;\Int)\cong{\rm Hom}(\thh_{k-1}(K_A;\Int),\Int)\cong \Int.$$
Next assume that $c_A \neq 0$. By (\ref{cohomk}) and (\ref{uct})
$$
0 \leftarrow {\rm Hom}(\hh_k(K_A;\Int),\Int) \leftarrow \Int^{|A|-1} \oplus \Int/d_A \Int \leftarrow {\rm Ext}(\hh_{k-1}(K_A;\Int),\Int) \leftarrow 0~.
$$
Therefore $\hh_k(K_A;\Int)\cong \Int^{|A|-1}$ and  ${\rm Ext}(\hh_{k-1}(K_A;\Int),\Int)= \Int/d_A \Int$.
Together with (\ref{EPH}) this implies that $\rk~ \thh_{k-1}(K_A;\Int)=0$. Therefore
$\thh_{k-1}(K_A;\Int)= \Int/d_A \Int$ and $\thh^{k-1}(K_A;\Int)=0$.
{\enp}
\noindent
{\bf Remark:} In the proof of (ii) it was observed that the function  $f \in \Int[\Int_n]$  given by
$f(\ell)=c_{\ell}$ if $0 \leq \ell \leq \varphi(n)$ and zero otherwise, is the image under $\theta^*$ of a  $k$-coboundary of $X(G)$.
This fact also appears (with a different proof) in Proposition 24 of \cite{MR10} and is attributed there to D. Fuchs.
\ \\ \\
{\bf Acknowledgment}
\\ I would like
to thank Vic Reiner for helpful discussions and comments.


\begin{thebibliography}{99}

\bibitem{Lang}
S. Lang, {\it Algebra}, Springer-Verlag, New York, 2002.

\bibitem{Munkres}
J.R. Munkres, {\it Elements of algebraic topology}, Addison-Wesley Publishing Company, Menlo
Park, CA, 1984.

\bibitem{MR10}
G. Musiker and V. Reiner, The cyclotomic polynomial topologically,
arXiv:1012.1844 .

\bibitem{Stanley}
R. Stanley, {\it Combinatorics and Commutative Algebra}, 2nd Edition, Birkh\"{a}user, Boston 1996.

\end{thebibliography}
\end{document}